\documentclass[reqno]{amsart}
\usepackage[centertags]{amsmath}
\usepackage{amsfonts}
\usepackage{amssymb}
\usepackage{amsthm}
\usepackage{newlfont}
\usepackage[top=0.9 in, bottom=0.9 in, left=0.9 in, right=0.9 in]{geometry}
\usepackage{mathrsfs}
\usepackage{dsfont}

\theoremstyle{plain}
\newtheorem{theo}{Theorem}[section]
\newtheorem{lemma}[theo]{Lemma}
\newtheorem{ax}{Axiom}
\newtheorem*{ax*}{Axiom}
\newtheorem*{conj}{Conjecture}
\newtheorem{cor}[theo]{Corollary}
\theoremstyle{remark}
\newtheorem{remark}[theo]{Remark}

\newcommand{\N}{\mathbb{N}}
\newcommand{\Z}{\mathbb{Z}}
\newcommand{\C}{\mathbb{C}}

\DeclareMathOperator{\Fix}{Fix}

\newcommand{\h}[2]{\mathscr{H}_{#1,#2}}
\newcommand{\q}[1]{\mathscr{Q}_{#1}}
\newcommand{\Q}{\mathscr{Q}}

\newcommand{\lb}{l}

\newcommand{\1}{\mathds{1}}

\newcommand{\set}[1]{\left\{ #1 \right\}}

\newcommand{\br}[1]{\left\langle #1 \right\rangle}

\newcommand{\parlengths}{\setlength{\parindent}{0pt}}
\begin{document}
\date{\today}

\title{FJRW-Rings and Landau-Ginzburg Mirror Symmetry in Two Dimensions}
\author{Pedro Acosta}
\address{Department of Mathematics, Brigham Young University, Provo, UT 84602, USA}
\email{pacosta@byu.net}

\maketitle\parlengths

\begin{abstract}
For any non-degenerate, quasi-homogeneous hypersurface singularity
$W$ and an admissible group of diagonal symmetries $G$, Fan,
Jarvis, and Ruan have constructed a cohomological field theory
which is a candidate for the mathematical structure behind the
Landau-Ginzburg A-model. When using the orbifold Milnor ring of a
singularity $W$ as a B-model, and the Frobenius algebra $\h WG$
constructed by Fan, Jarvis, and Ruan, as an A-model, the following
conjecture is obtained: For a quasi-homogeneous singularity $W$
and a group $G$ of symmetries of $W$, there is a dual singularity
$W^T$ such that the orbifold A-model of $W/G$ is isomorphic to the
B-model of $W^T$. I will show that this conjecture holds for a
two-dimensional invertible loop potential $W$ with its maximal
group of diagonal symmetries $G_{W}$.

\end{abstract}

\section{Introduction}\label{intro}

In a recent paper \cite{FJR}, Fan, Jarvis, and Ruan constructed
the mathematical theory (FJRW-theory) behind the Landau-Ginzburg
A-model. Their construction gives, among other things, a Frobenius
algebra $\h WG$ which is determined by a non-degenerate
quasi-homogeneous polynomial $W$ and an admissible group $G$ of
diagonal symmetries of $W$.

Landau-Ginzburg models have been studied extensively in the
physics literature (see for example \cite{green}). In particular,
a mirror construction for these models was suggested by Berglund
and H\"{u}bsch \cite{berglund}. This mirror construction used the
so-called `invertible' potentials: quasi-homogeneous polynomials
with the same number of monomials and variables. In the FJRW
theory, these potentials must satisfy non-degeneracy conditions
given by,
\begin{enumerate}
    \item The potential $W$ must have an isolated singularity at
    the origin.
    \item The weights (charges) of $W$ must be uniquely
    determined.
\end{enumerate}

These non-degeneracy conditions imply that the invertible
potentials will be of two kinds \cite{kreuzer2},
\begin{equation}
W_{loop}=X_{1}^{a_{1}}X_{2}+\dots+X_{N-1}^{a_{N-1}}X_{N}+X_{N}^{a_{N}}X_{1}
\end{equation}
\begin{equation}
W_{chain}=X_{1}^{a_{1}}X_{2}+\dots+X_{N-1}^{a_{N-1}}X_{N}+X_{N}^{a_{N}}
\end{equation}

In \cite{kreuzer1}, Kreuzer proved that the Berglund-H\"{u}bsch
construction can be used to show that the types of potentials
previously described satisfy a certain type of mirror symmetry.
Based on Kreuzer's work, Krawitz \cite{krawitz} conjectured that
the FJRW construction satisfies the following,

\begin{conj}{(Landau-Ginzburg Mirror Symmetry Conjecture):}
For a non-degenerate, quasi-homogeneous, invertible singularity
$W$ and its maximal group of diagonal symmetries $G_{W}$, there is
a dual singularity $W^T$ such that the FJRW-ring of
$\mathscr{H}_{W,G_{W}}$ is isomorphic to the (unorbifolded) Milnor
ring of $W^T$.
\end{conj}

This conjecture has been verified for simple singularities by Fan,
Jarvis, and Ruan in \cite{FJR}, and for unimodal and bimodal
singularities by Priddis et al. in \cite{priddis}. Recently, Fan
and Shen \cite{fan} proved that this conjecture is true for all
chain-type potentials in two dimensions $(N=2)$. Following some of
the ideas of Kreuzer, I will show that the Landau-Ginzburg Mirror
Symmetry conjecture is true for loop-type potentials in two
dimensions. Combined with the work of Fan and Shen, this completes
the proof of the conjecture for all superpotentials in two
dimensions. In a recent preprint \cite{krawitz}, Krawitz has shown
us a proof of the conjecture for $N=3$ and higher.

\subsection{Outline of the Paper}\label{outline}
The organization of the paper is as follows. First, a review of
the FJRW construction will be given in Section \ref{review}.
Section \ref{additional} describes some additional notation that
will be used throughout the paper. In Section \ref{loops}, a proof
of the conjecture will be given for loop potentials in two
dimensions.

\subsection{Acknowledgments}\label{acknow}
The author would like to thank T. Jarvis for helpful discussions and
comments. He also thanks M. Krawitz and N. Priddis for fruitful
discussions. Special acknowledgments are given to the Department of
Mathematics of Brigham Young University for their support.

\subsection{Review of the Construction}\label{review}
In this section, I will present a simplified description of the
FJRW-theory, similar to the one outlined in \cite{priddis}. For a
full description of the theory, the reader is referred  to the
original paper of Fan, Jarvis, and Ruan \cite{FJR}.

A \emph{quasi-homogeneous} polynomial $W \in
\C[X_{1},\dots,X_{N}]$ is defined as a polynomial for which there
exist rational positive degrees $q_{1},\dots,q_{N} \in
\mathbb{Q}^{>0}$, such that for any $\lambda \in \C^{*}$
\[
W(\lambda^{q_{1}}X_{1},\dots,\lambda^{q_{N}}X_{N})=W(X_{1},\dots,X_{N}).
\]

For $i \in \{1,...,N\}$, we will call $q_{i}$ the \emph{weight} of
the variable $X_{i}$.

Let $W: \C^N\longrightarrow\C$ be a quasi-homogeneous polynomial
satisfying the non-degeneracy conditions (1) and (2). We define
the \emph{local algebra of $W$}, also known as the \emph{Milnor
ring}, by
\begin{equation}
\Q_{W}:=\C[X_{1},\dots,X_{N}]/\text{Jac}(W),
\end{equation}
where Jac(W) is the Jacobian ideal of W, generated by partial
derivatives
\begin{equation}
\text{Jac}(W):=\left(\frac{\partial{W}}{\partial{X_{1}}},\dots,\frac{\partial{W}}{\partial{X_{N}}}\right)
\end{equation}
It is easy to see that the local algebra is generated by monomials
of the form $\prod X_{i}^{b_{i}}$. The local algebra is graded by
the weighted-degree of each monomial, where $X_{i}$ has weight
$q_{i}$.

The local algebra contains a unique highest-degree element given
by
$\text{det}\left(\frac{\partial^2{W}}{\partial{X_{i}}\partial{X_{j}}}\right)$,
whose degree is given by
\begin{equation}
\hat{c}_W=\sum_{i=1}^{N}(1-2q_{i}).
\end{equation}

Also, the dimension of the local algebra as a vector space over
$\C$ is given by
\begin{equation}
\mu=\prod_{i=1}^{N} \left(\frac{1}{q_{i}}-1\right).
\end{equation}

For $f,g \in \Q_{W}$, we define the \emph{residue pairing}
$\br{f,g}:\mathscr{Q}_W\times \mathscr{Q}_W \longrightarrow \C $
by
\[
fg =
\frac{\br{f,g}}{\mu}\text{det}\left(\frac{\partial^2{W}}{\partial{X_{i}}\partial{X_{j}}}\right)
+ \text{ lower order terms.}
\]
This pairing endows the local algebra with the structure of a
Frobenius algebra. This structure is known as the
\emph{unorbifolded Landau-Ginzburg B-model}.

In order to define the FJRW-ring $\h WG$, or A-model, we need both
a potential $W$ and an admissible group of diagonal symmetries $G$
of $W$. For the definition of admissible, see \cite{FJR}. We
define the \emph{maximal group $G_{W}$ of diagonal symmetries} of
W by
\begin{equation}
G_{W}:=\{(\alpha_{1},...,\alpha_{N})\in (\C^{*})^{N} \,|\,
W(\alpha_{1}X_{1},...,\alpha_{N}X_{N})=W(X_{1},...X_{N})\}
\end{equation}
This group is known to be admissible and, as shown in \cite{FJR},
when $W$ satisfies the non-degeneracy conditions, the group
$G_{W}$ is finite. Also note that the element $J$, defined by
\begin{equation}
J:=\left(e^{2\pi i q_{1}},\dots,e^{2\pi i q_{N}}\right)
\end{equation}
always belongs to $G_{W}$, and the group $\br{J}$ generated by $J$
is always admissible.

The original definition of the state space of $\h WG$ is in terms
of Lefschetz thimbles. However, it can be described more simply in
terms of sums of local algebras.

For $g \in G$, let $\Fix$ g $\subseteq \C^N$ be the set of fixed
points of g, and let $N_{g}$ be its dimension. Also, let $W |\Fix
g$ be the potential restricted to the fixed point locus of $g$.
Define the vector space $\mathscr{H}_g$ by
\begin{equation}
\mathscr{H}_g:=\Q_{W |\Fix g} \cdot \omega,
\end{equation}
where $\omega=dX_{i_{1}}\wedge dX_{i_{2}}\wedge\dots\wedge
dX_{i_{N_{g}}}$.

The \emph{state space} of $\mathscr{H}_{W,G}$ is defined as the
G-invariant subspace of the sum of the $\mathscr{H}_{g}$:
\begin{equation}
\mathscr{H}_{W,G} := \left(\bigoplus_{g\in
G}\mathscr{H}_g\right)^{G}.
\end{equation}

This space can be graded by the so-called \emph{W-degree}. In
order to define the W-degree of each element in
$\mathscr{H}_{W,G}$, we note that any element $g \in G$ can be
written in the form
\[
g=\left(e^{2 \pi i\theta_{1}^g},\dots,e^{2 \pi
i\theta_{N}^g}\right).
\]

If the phases $\theta_{j}^{g}$ satisfy the condition $0 \leq
\theta_{j}^g < 1$, then we denote them by $\Theta_{j}^{g}$. Note
that any element $g \in G$ can be uniquely written in the form
\begin{equation}
g=\left(e^{2 \pi i\Theta_{1}^g},\dots,e^{2 \pi
i\Theta_{N}^g}\right),\quad \text{with}\quad 0\leq
\Theta_{j}^{g}<1.
\end{equation}

Let $\alpha_g \in \left(\mathscr{H}_g\right)^G$, then we define
the W-degree of $\alpha_g$ by
\begin{equation}
\deg_W(\alpha_g):=N_g +2\sum_{j=1}^N (\Theta_{j}^g-q_{j})
\end{equation}

From (11), it is easy to check that $\Fix g=\Fix g^{-1}$ and thus
there is a canonical isomorphism $I:\mathscr{H}_g \rightarrow
\mathscr{H}_{g^{-1}}$. From this we can see that the pairing on
$\q{W|_{\Fix g}}$ induces a pairing $\eta_{g}$
\[
\eta_{g}:(\mathscr{H}_g)^G \otimes
(\mathscr{H}_{g^{-1}})^G\to\C,\quad \text{given by}\quad
\eta_{g}(a,b)=\br{a,I^{-1}(b)}.
\]
The pairing on $\h WG$ is  the direct sum of the pairings
$\eta_{g}$. Fixing a basis for $\h WG$, we denote the pairing by a
matrix $\eta_{\alpha,\beta}=\br{\alpha,\beta}$, with inverse
$\eta^{\alpha,\beta}$.

For each pair of non-negative integers $g$ and $n$, with
$2g-2+n>0$, the Fan-Jarvis-Ruan-Witten (FJRW) cohomological field
theory produces classes
$\Lambda_{g,n}^W(\alpha_1,\alpha_2,\dots,\alpha_n)\in
H^*(\overline {\mathscr{M}}_{g,n})$  of complex codimension $D$
for each $n$-tuple $(\alpha_1,\alpha_2,\dots,\alpha_n)\in (\h
WG)^n$. Here, $\overline{\mathscr{M}}_{g,n}$ is the stack of
stable curves of genus $g$ with $n$ marked points, and the
codimension $D$ is given by
\[
D:=\hat c_W(g-1)+\frac 12\sum_{i=1}^n\deg_W(\alpha_i),
\]
We define \emph{three-point correlators} by
\[
\br{\alpha_1,\alpha_2,\alpha_3}:=\int_{\overline
{\mathscr{M}}_{0,3}}\Lambda_{0,3}^W(\alpha_1,\alpha_2,\alpha_3).
\]
The three-point correlator $\br{\alpha_1,\alpha_2,\alpha_3}$
vanishes unless $D$ is zero. These three-point correlators can be
used to define structure constants for a multiplication on
$\mathscr{H}_{W,G_{W}}$. If $r,s\in \h{W}{G}$, this multiplication
is defined by
\begin{equation}
r\star
s:=\sum_{\alpha,\beta}\br{r,s,\alpha}\eta^{\alpha,\beta}\beta,
\label{mult}
\end{equation}
where the sum is taken over all choices of $\alpha$ and $\beta$ in
a fixed basis of $\h WG$.

As described in \cite{priddis}, in genus zero with three marked
points, the class $\Lambda_{0,3}^W(\alpha_1,\alpha_2,\alpha_3)$
satisfies the following axioms that allow us to compute most of
the three-point correlators $\br{\alpha_1,\alpha_2,\alpha_3}$
explicitly.

\begin{ax} Dimension: If $2D\notin  \Z$, then $\Lambda_{g,n}^W(\alpha_1,\alpha_2,\dots,\alpha_n)=0$. Otherwise, $2D$ is the real codimension of the class $\Lambda_{g,n}^W(\alpha_1,\alpha_2,\dots,\alpha_n)$. In particular, if $g=0$ and $n=3$, then $\br{\alpha_1,\alpha_2,\alpha_3}=0$ unless $D=0$.
\end{ax}

Notice that in the case where $g=0$ and $n=3$, if $D$ is the
codimension of the class
$\Lambda_{0,3}^W(\alpha_1,\alpha_2,\alpha_3)$, then $D=0$ if and
only if $\sum_{i=1}^3\deg_W\alpha_i=2\hat c$.

\begin{ax} Symmetry:
Let $\sigma\in S_n$. Then
\[
\Lambda_{g,n}^W(\alpha_1,\alpha_2,\dots,\alpha_n)=\Lambda_{g,n}^W(\alpha_{\sigma(1)},\alpha_{\sigma(2)},\dots,\alpha_{\sigma(n)})
\]
\end{ax}

The next few axioms rely on the degrees of line bundles
$\mathscr{L}_1,\dotsc,\mathscr{L}_N$ endowing an orbicurve with a
so-called \emph{$W$-structure}; however, this can be reduced to a
simple numerical criterion. Consider the class
$\Lambda_{g,n}^W(\alpha_1,\alpha_2,\dots,\alpha_k)$, with
$\alpha_j\in(\mathscr{H}_{h_j})^G$ for each $j$. For each variable
$X_j$, define $\lb_j$ to be the \emph{degree of the line bundle}
$\mathscr{L}_{j}$. By \cite{FJR}, this degree is given by the
equation
\[
\lb_j=q_j(2g-2+k)-\sum_{i=1}^k\Theta_j^{h_i}.
\]

\begin{ax} Integer degrees: If $\lb_j\notin \Z$ for some $j\in\set{1,\dots,N}$, then $\Lambda_{g,n}^W(\alpha_1,\alpha_2,\dots,\alpha_n)=~0$.
\end{ax}

\begin{ax}
Concavity: If $\lb_{j}<0$ for all $j\in\set{1,2,3}$, then
$\br{\alpha_1,\alpha_2,\alpha_3}=1$.
\end{ax}

The next axiom is related to the Witten map:
\begin{eqnarray*}
\mathcal W:\bigoplus_{\substack{j=1\\}}^N H^{0}(\mathscr{C},\mathscr{L}_{j})\rightarrow \bigoplus_{j=1}^N H^{1}(\mathscr{C},\mathscr{L}_{j})\\
\mathcal W=\left(\overline{\frac{{\partial W}}{{\partial x_1}}},
\overline{\frac{{\partial W}}{{\partial x_2}}}, \dots,
\overline{\frac{{\partial W}}{{\partial x_N}}}\right)
\end{eqnarray*}

The dimensions of the cohomologies
$H^{0}(\mathscr{C},\mathscr{L}_{j})$ and
$H^{1}(\mathscr{C},\mathscr{L}_{j})$ are $h^0_j$ and $h^1_j$
respectively. These dimensions are known to be given by
\[
h^0_j:=\begin{cases}
0       & \text{if } \lb_j< 0\\
\lb_j+1 & \text{if }\lb_j \geq 0
\end{cases}
\]
\[
h^1_j:=\begin{cases}
-\lb_j-1 & \text{if } \lb_j<0\\
0        & \text{if } \lb_j\geq 0
\end{cases}
\]
so that both are non-negative integers satisfying $h^0_j - h^1_j =
l_j + 1$. Moreover, we have that $D=\sum_{j=1}^N (h^0_j - h^1_j)$.

In $\overline {\mathscr{M}}_{g,n}$, if
$\Lambda_{g,n}^W(\alpha_1,\alpha_2,\dots,\alpha_n)$ is a class of
codimension zero, then these classes are constant and so, abusing
notation, we will simply consider
$\Lambda_{g,n}^W(\alpha_1,\alpha_2,\dots,\alpha_n)$ to be a
complex number. We will use this convention through the rest of
the paper.

\begin{ax}
Index Zero: Consider the class
$\Lambda_{g,n}^W(\alpha_1,\alpha_2,\dots,\alpha_n)$, with
$\alpha_i\in \h{\gamma_i}{G}$. If $\Fix \gamma_i=\set 0$ for each
$i\in\set{1,2,\dots,n}$ and
\[
D=\sum_{j=1}^N (h^0_j - h^1_j) = 0,
\]
then $\Lambda_{g,n}(\alpha_1,\alpha_2,\dots,\alpha_n)$ is of
codimension zero and
$\Lambda_{g,n}^W(\alpha_1,\alpha_2,\dots,\alpha_n)$ is equal to
the degree of the Witten map.
\end{ax}

\begin{ax}
Composition: If the four-point class,
$\Lambda_{g,n}^W(\alpha_1,\alpha_2,\alpha_3,\alpha_4)$ is of
codimension zero, then it decomposes as sums of three-point
correlators in the following way:
\[
\Lambda_{0,4}^W(\alpha_1,\alpha_2,\alpha_3,\alpha_4)=\sum_{\beta,\delta}\br{\alpha_1,
\alpha_2,\beta}\eta^{\beta,\delta}\br{\delta,\alpha_3,\alpha_4}=\sum_{\beta,\delta}\br{\alpha_1,
\alpha_3,\beta}\eta^{\beta,\delta}\br{\delta,\alpha_2,\alpha_4}.
\]
\end{ax}

Note that $\Fix J=\set 0$ so $\mathscr{H}_{J}\cong\C$. Let $\1$ be
the element in $\mathscr{H}_{J}$ corresponding to $1\in\C$. This
element has $\deg_W(\1)=0$ and it turns out to be the identity
element in the FJRW-ring. The next axiom deals with this element.

\begin{ax}
Pairing: For any $\alpha_1,\alpha_2\in \h WG$, we have
$\br{\alpha_1,\alpha_2,\mathds{1}}=\eta_{\alpha_1,\alpha_2}$.
\end{ax}

\begin{ax} \label{ax:sums} Sums of singularities: If $W_1\in\C[X_1,\dots,X_r]$ and $W_2\in\C[Y_1,\dots,Y_s]$ are two non-degenerate, quasi-homogeneous polynomials with maximal symmetry groups $G_1$ and $G_2$, then the maximal symmetry group of $W=W_1+W_2$ is $G=G_1\times G_2$, and there is an isomorphism of Frobenius algebras
\[
\h{W}{G}\cong \h{W_1}{G_{W_1}}\otimes \h{W_2}{G_{W_2}}
\]
\end{ax}

\subsection{Additional Notation}\label{additional}
Throughout this paper we adopt the following notation. Let $g \in
G_{W}$. If $\Fix g=\{0\}$, define
\[
\textbf{e}_g:=1\in \mathscr{H}_{g} \cong \C,
\]
otherwise, if $g$ fixes the variables $ X_{i_1},\dots,
X_{i_{N_g}}$, define
\[
\textbf{e}_{g}:=dX_{i_1}\wedge dX_{i_2}\wedge\dots\wedge
dX_{i_{N_g}}\in \mathscr{H}_{g}.
\]

The identity element of $G_{W}$ will be denoted by $i_{W}$.

\section{Two-Dimensional Loop Potentials}\label{loops}

In this section we show that the Landau-Ginzburg Mirror Symmetry
Conjecture holds for the so-called loop potentials in two
dimensions. These potentials are of the form
\begin{equation}
W=X_{1}^{a_{1}}X_{2}+X_{2}^{a_{2}}X_{1},
\end{equation}

where $a_{1}$, $a_{2}$ $\in \N^{>1}$. Note that quasi-homogeneity
implies that $a_{1}q_{1}+q_{2}=1$ and $a_{2}q_{2}+q_{1}=1$ for the
weights $q_{1}$, $q_{2}$ of the loop potential, so
$q_{1}=\frac{a_{2}-1}{a_{1}a_{2}-1}$ and
$q_{2}=\frac{a_{1}-1}{a_{1}a_{2}-1}$.

\begin{theo}
For an arbitrary two-dimensional loop potential
\[W=X_{1}^{a_{1}}X_{2}+X_{2}^{a_{2}}X_{1}\] the
group of diagonal symmetries $G_{W}$ is cyclic of order
\[
|G_{W}|=a_{1}a_{2}-1.
\]
Moreover, $G_{W}$ can be generated by either one of the following
elements: $g_{1}:=(e^{2\pi i \theta_{1}^{(1)}},e^{2\pi i
\theta_{2}^{(1)}})$ or \\
 $g_{2}:=(e^{2\pi i\theta_{1}^{(2)}},e^{2\pi i \theta_{2}^{(2)}})$,
 where
\begin{equation}
\theta_{1}^{(1)}:=\frac{-1}{|G_{W}|}, \quad
\theta_{2}^{(1)}:=\frac{a_{1}}{|G_{W}|}, \quad
\theta_{1}^{(2)}:=\frac{a_{2}}{|G_{W}|}, \quad
\theta_{2}^{(2)}:=\frac{-1}{|G_{W}|}.
\end{equation}
\end{theo}

\begin{proof}

From the definition of the maximal group of diagonal symmetries,
we see that any $g \in G_{W}$ can be expressed in the form
$g=(\alpha_{1},\alpha_{2})$, satisfying the conditions
\begin{equation}
\alpha_{1}^{a_{1}}\alpha_{2}=1\quad \text{and}\quad
\alpha_{2}^{a_{2}}\alpha_{1}=1.
\end{equation}

This implies that $\alpha_{1}^{a_{1}a_{2}-1}=1$ and
$\alpha_{2}^{a_{1}a_{2}-1}=1$. Thus, $\alpha_{1}$ and $\alpha_{2}$
are primitive roots of unity of order $a_{1}a_{2}-1$. From (16) it
is clear that determining the value of $\alpha_{1}$ fixes the
value of $\alpha_{2}$, and therefore $G_{W}$ must be isomorphic to
the additive group of integers modulo $a_{1}a_{2}-1$. Hence,
$G_{W}$ is cyclic of order $|G_{W}|=a_{1}a_{2}-1$.

From (15) it is clear that $g_{2}$, $g_{2} \in G_{W}$, and that
they are generators of the group.

\end{proof}

\begin{remark}
Note that $g_{1}^{-a_{2}}=g_{2}$ and $g_{2}^{-a_{1}}=g_{1}$.
\end{remark}

We now want to construct the state space of $\h GW$. As mentioned
in Section \ref{review} it is possible to do this in terms of
Milnor rings.

\begin{theo}
The state space of the FJRW-ring of a two-dimensional loop
potential W with maximal group diagonal of symmetries $G_{W}$ is
given by
\[
\mathscr{H}_{W,G_{W}}=\left(X_{1}^{a_{1}-1}\right)\textbf{e}_{i_{W}}\oplus
\left(X_{2}^{a_{2}-1}\right)\textbf{e}_{i_{W}}
\bigoplus_{\tiny\begin{array}{c}g \in G_{W}\\ g \neq
i_{W}\end{array}}\C \textbf{e}_{g}
\]
\end{theo}

\begin{proof}

Since $G_{W}=\br{g_{i}}$ for $i \in \{1,2\}$, we have that we can
construct the state space of the FJRW-ring by taking powers of any
of the generators $g_{i}$, i.e.
\[
\mathscr{H}_{W,G_{W}} =
\left(\bigoplus_{k=0}^{|G_{W}|-1}\mathscr{H}_{g_i^{k}}\right)^{G_{W}}.
\]

By looking at (15), it is not hard to see that $|G_{W}|$ and
$|G_{W}|\theta_{j}^{(i)}$ are relatively prime for any of the
phases $\theta_{j}^{(i)}$ of $g_{i}$. From this result we see that
$\Fix g_{i}^{k}=\{0\}$ for $k \neq 0$, and therefore
$N_{g_{i}^k}=0$. Thus, we only need to find the elements in
$\mathscr{H}_{g_{i}^{0}}$ that are invariant under the action of
$G_{W}$. We first note that $\Fix {g_{i}^{0}}=\Fix {i_{W}}=\C^2$,
where $i_{W}$ is the identity of $G_{W}$. In order to find
$\mathscr{H}_{i_{W}}$ then, we need to compute the Milnor ring of
$W$. We recall that this ring is given by,
\[
\Q_{W}=\C[X_{1},X_{2}]/\text{Jac}(W)
\]

where $\text{Jac}(W)$ is the Jacobian ideal of $W$. From these
relations, we find that a basis for $\Q_{W}$ as a vector space
over $\C$ is given by monomials of the form
\begin{equation}
X_{1}^{b_{1}}X_{2}^{b_{2}}
\end{equation}
where $0 \leq b_{1} < a_{1}$ and $0 \leq b_{2} < a_{2}$.

As described in Section \ref{review}, a basis for
$\mathscr{H}_{i_{W}}$ is therefore given by elements of the form
\[
\left( X_{1}^{b_{1}}X_{2}^{b_{2}}\right)\textbf{e}_{i_{W}},\quad 0
\leq b_{1} < a_{1}, \quad 0 \leq b_{2} < a_{2}.
\]

The elements in $\mathscr{H}_{i_{W}}$ invariant under the action
of $G_{W}$ must satisfy
\begin{equation}
\sum_{j=1}^{2} \theta_{j}^{(i)}b_{j} + \sum_{j=1}^2
\theta_{j}^{(i)} = m ,\,  \text{where} \, m \in \Z
\end{equation}

This relation must hold true for any generator $g_{i}$ of $G_{W}$,
so in particular we can pick $i=1$ in the above condition. Since
$\theta_{1}^{(1)}$ is negative and $\theta_{2}^{(1)}$ is positive,
we have that the maximum and minimum values of (20) will be
attained  by $\left(X_{2}^{a_{2}-1}\right)\textbf{e}_{i_{W}}$ and
$\left(X_{1}^{a_{1}-1}\right)\textbf{e}_{i_{W}}$, respectively.
The value of m in these cases will be $1$ and $0$ respectively. If
we take any other element of $\mathscr{H}_{i_{W}}$ the value of
(18) will be strictly between $0$ and $1$. Therefore, the only
elements of $\mathscr{H}_{i_{W}}$ fixed under the action of
$G_{W}$ are $\left(X_{1}^{a_{1}-1}\right)\textbf{e}_{i_{W}} $ and
$\left( X_{2}^{a_{2}-1}\right)\textbf{e}_{i_{W}}$.

\end{proof}

\begin{remark}
The dimension of $\mathscr{H}_{W,G_{W}}$ as a vector space over
$\C$ is given by $a_{1}a_{2}$.
\end{remark}

Now that we have constructed the state space for
$\mathscr{H}_{W,G_{W}}$, we would like to find the potential
$W^{T}$ that will be the mirror dual of $W$. From the
Berglund-H\"{u}bsch mirror construction \cite{berglund}, this dual
potential is given by

\begin{equation}
W^{T}=\overline{X}_{1}\overline{X}_{2}^{a_{1}}+\overline{X}_{2}\overline{X}_{1}^{a_{2}}.
\end{equation}

Let $\overline{q}_{i}$ be the weight of the variable
$\overline{X}_{i}$ in $W^{T}$. It is not hard to check that
\begin{equation}
a_{1}\overline{q}_{2}+\overline{q}_{1}=1 \quad\text{and}\quad
a_{2}\overline{q}_{1}+\overline{q}_{2}=1,
\end{equation}
and that
\begin{equation}
\overline{q}_{i}:=\theta_{1}^{(i)}+\theta_{2}^{(i)}, \quad i \in
\{1,2\}.
\end{equation}

 We claim that this new potential $W^{T}$ is the mirror
dual of $W$, which means that $\left(\mathscr{H}_{W,G_{W}},\star
 \, \right)\cong \mathscr{Q}_{W^{T}}$ and $\left(\mathscr{H}_{W^{T},G_{W^{T}}},\star
 \, \right)\cong \mathscr{Q}_{W}$, where $\star$ is the multiplication defined in (13).

\begin{theo}
For a loop potential $W=X_{1}^{a_{1}}X_{2}+X_{1}X_{2}^{a_{2}}$
with maximal group of diagonal symmetries $G_{W}$, there is a dual
potential $W^{T}$ given by
\[
W^{T}=\overline{X}_{1}\overline{X}_{2}^{a_{1}}+\overline{X}_{2}\overline{X}_{1}^{a_{2}}
\]
such that $\mathscr{H}_{W,G_{W}}\cong \mathscr{Q}_{W^{T}}$ as
graded Frobenius algebras, where $\mathscr{H}_{W,G_{W}}$ is graded
by W-degree and $\mathscr{Q}_{W^{T}}$ is graded by the weighted
degree of monomials.
\end{theo}

In order to prove Theorem 2.5, we must first prove a series of
results:

\begin{lemma}
Every element $g \in G_{W}$ can be written in the form
$g=Jg_{1}^{\alpha}g_{2}^{\beta}$ in a unique way, where $0 \leq
\alpha \leq a_{2}-1$ and $0 \leq \beta \leq a_{1}-1$, with the
exception of
 $i_{W}=Jg_{1}^{a_{2}-1}=Jg_{2}^{a_{1}-1}$.
\end{lemma}

\begin{proof}

First, note that
\[
Jg_{1}^{a_{2}-1}=\left(e^{2\pi i
(q_{1}+(a_{2}-1)\theta_{1}^{(1)})},e^{2\pi i
(q_{2}+(a_{2}-1)\theta_{2}^{(1)})} \right)=i_{W},\quad  \text{and
that}
\]
\[
Jg_{2}^{a_{1}-1}=\left(e^{2\pi i
(q_{1}+(a_{1}-1)\theta_{1}^{(2)})},e^{2\pi i
(q_{2}+(a_{1}-1)\theta_{2}^{(2)})} \right)=i_{W}.
\]

Now suppose that for some $g \in G_{W}$ with $g \neq i_{W}$, we
have that
$g=Jg_{1}^{\alpha_{1}}g_{2}^{\beta_{1}}=Jg_{1}^{\alpha_{2}}g_{2}^{\beta_{2}}$,
where $0 \leq \alpha_{1}, \alpha_{2} \leq a_{2}-1$ and $0 \leq
\beta_{1}, \beta_{2} \leq a_{1}-1$. Assume without loss of
generality that $\alpha_{2}\geq \alpha_{1}$. If we divide one
representation by the other, we find that
\[
i_{W}=g_{1}^{\alpha_{2}-\alpha_{1}}g_{2}^{\beta_{2}-\beta_{1}}=g_{1}^{\alpha_{2}-\alpha_{1}}g_{1}^{-a_{2}(\beta_{2}-\beta_{1})},
\]
where the last equality comes after invoking \emph{Remark 2.2}.
This implies that
$\alpha_{2}-\alpha_{1}-a_{2}(\beta_{2}-\beta_{1})=m|G_{W}|$, where
$m$ is an integer. It is not hard to show that
\[
-a_{2}(a_{1}-1) \leq
\alpha_{2}-\alpha_{1}-a_{2}(\beta_{2}-\beta_{1}) \leq a_{2}-1
+a_{2}(a_{1}-1)=|G_{W}|,
\]
and therefore, the only possible values that $m$ can take are 0
and 1. The only way in which $m=1$ is by letting
$\alpha_{2}=a_{2}-1$, $\alpha_{1}=0$, $\beta_{2}=0$, and
$\beta_{1}=a_{1}-1$. However, this would mean that $g=i_{W}$,
which is impossible.

It is straightforward to show that the only way in which $m=0$ is
that $\alpha_{1}=\alpha_{2}$ and that $\beta_{1}=\beta_{2}$, but
this means that the representation of $g$ in the form
$Jg_{1}^{\alpha}g_{2}^{\beta}$ is unique.

We have thus far shown that there are exactly $a_{1}a_{2}-1$
different elements of $G_{W}$ that can be written in the form
$Jg_{1}^{\alpha}g_{2}^{\beta}$, with $0 \leq \alpha \leq a_{2}-1$
and $0 \leq \beta \leq a_{1}-1$. Since the order of $G_{W}$ is
also $a_{1}a_{2}-1$, then every element of $G_{W}$ can be written
uniquely in the form $Jg_{1}^{\alpha}g_{2}^{\beta}$, with the
exception of $i_{W}.$
\end{proof}

\begin{cor}
Let $\gamma \in \mathscr{H}_{W,G_{W}}$, and suppose that $\gamma
\in \mathscr{H}_{g}$ for some $g \in G_{W}$, where
$g=Jg_{1}^{\alpha}g_{2}^{\beta}$, $0 \leq \alpha \leq a_{2}-1$, $0
\leq \beta \leq a_{1}-1$. Then the W-degree of $\gamma$ is given
by
\[
\deg_{W}(\gamma)=2(\alpha\overline{q}_{1}+\beta\overline{q}_{2}).
\]
\end{cor}
\begin{proof}
We divide the proof in two cases: $g=i_{W}$ and $g \neq i_{W}$.

\emph{Case 1}: Suppose that $\gamma \in \mathscr{H}_{i_{W}}$.
Then, the W-degree of $\gamma$ will be given by

\[
\deg_{W}(\gamma)=N_{i_{W}}+2(0-q_{1})+2(0-q_{2})
\]
\[
=2-2q_{1}-2q_{2}.
\]

A simple computation shows that
$\deg_W(\gamma)=2(a_{2}-1)\overline{q}_{1}=2(a_{1}-1)\overline{q}_{2}.$

\emph{Case 2}: Suppose that $g \neq i_{W}$ and that $\gamma \in
\mathscr{H}_{g}$. Let $g=Jg_{1}^{\alpha}g_{2}^{\beta}$, where $0
\leq \alpha \leq a_{2}-1$, $0 \leq \beta \leq a_{1}-1$. Then
\[
\theta_{1}^{g}=q_{1}+\alpha\theta_{1}^{(1)}+\beta\theta_{1}^{(2)}
\quad \text{and} \quad
\theta_{2}^{g}=q_{2}+\alpha\theta_{2}^{(1)}+\beta\theta_{2}^{(2)}.
\]

It is not hard to show that $0 \leq
q_{1}+\alpha\theta_{1}^{(1)}+\beta\theta_{1}^{(2)} \leq 1$ and
that $0 \leq q_{2}+\alpha\theta_{2}^{(1)}+\beta\theta_{2}^{(2)}
\leq 1$. Note that the only time that $\theta_1^{g}=1$ or
$\theta_{2}^{g}=1$ is when $g=i_{W}$, but this was considered in
\emph{Case 1}, so we will assume that $0 \leq
\theta_{1}^{g},\theta_{2}^{g} <1$. Therefore, we have that
$\Theta_{1}^{g}=\theta_{1}^{g}$ and
$\Theta_{2}^{g}=\theta_{1}^{g}$. We can now use (12) to compute
the W-degree of $\gamma$
\[
\deg_W(\gamma)=N_{\gamma}+2(q_{1}+\alpha\theta_{1}^{(1)}+\beta\theta_{1}^{(2)}-q_{1})+2(q_{2}+\alpha\theta_{2}^{(1)}+\beta\theta_{2}^{(2)}-q_{2})
\]
\[
=0+2\alpha(\theta_{1}^{(1)}+\theta_{2}^{(1)})+2\beta(\theta_{1}^{(2)}+\theta_{2}^{(2)})=2(\alpha\overline{q}_{1}+\beta\overline{q}_{2}).
\]
\end{proof}

\begin{lemma}
For any integer $c$ with  $0 \leq c < a_{1}-1$, we have that
$(\textbf{e}_{Jg_{i}})^c=\textbf{e}_{Jg_{i}^{c}}$, where $i \in
\{1,2\}$.

\end{lemma}
\begin{proof}

For $c=0$ the result is trivial since
\[
(\textbf{e}_{Jg_{i}})^0=\mathds{1}=\textbf{e}_J=\textbf{e}_{Jg_{i}^{0}}.
\]

Now suppose that for $1 \leq c < a_{1}-1$, we have that
$(\textbf{e}_{Jg_{i}})^{c-1}=\textbf{e}_{Jg_{i}^{c-1}}$, and
consider the product $(\textbf{e}_{Jg_{i}})^{c-1}\star\,  \,
\textbf{e}_{Jg_{i}}$. By definition, this product will be given by
$\sum_{\alpha,\beta}\br{(\textbf{e}_{Jg_{i}})^{c-1},\textbf{e}_{Jg_{i}},\alpha}\eta^{\alpha,\beta}\beta$.
Using our assumption, we find that
\begin{equation}
(\textbf{e}_{Jg_{i}})^{c-1}\star\, \,
\textbf{e}_{Jg_{i}}=\sum_{\alpha,\beta}\br{\textbf{e}_{Jg_{i}^{c-1}},\textbf{e}_{Jg_{i}},\alpha}\eta^{\alpha,\beta}\beta.
\end{equation}
For these correlators to be non-zero we need
$\deg_{W}(\textbf{e}_{jg_{i}^{c-1}})+\deg_{W}(\textbf{e}_{jg_{i}})+\deg_{W}(\alpha)=2\hat{c}$.
Using \emph{Corollary 2.7}, we find that this last relation is
equivalent to
\[
2c\overline{q}_{i}+\deg_W(\alpha)=2(a_{i+1}-1)\overline{q}_{i}+2(a_{i}-1)\overline{q}_{i+1}
\]
\[
\Rightarrow
\deg_W(\alpha)=2(a_{i+1}-1-c)\overline{q}_{i}+2(a_{i}-1)\overline{q}_{i+1}.
\]

From \emph{Corollary 2.7}, it is easy to see that if
$\gamma=Jg_{i}^{a_{i+1}-1-c}g_{i+2}^{a_{i}-1}$, then
$\deg_W(\textbf{e}_\gamma)=\deg_{W}(\alpha)$, and since $\gamma
\neq i_W$, there is only one basis element coming from
$\mathscr{H}_{\gamma}$, and thus, the sum in equation (22) reduces
to a single term. Also note that $\gamma g_{i}^{c}=i_{W}$, and
therefore, $\textbf{e}_{\gamma}$ pairs up with
$\textbf{e}_{Jg_{i}^c}$, which gives
\[
(\textbf{e}_{Jg_{2}})^{c}
=\br{\textbf{e}_{Jg_{i}^{c-1}},\textbf{e}_{Jg_{i}},\textbf{e}_{\gamma}}\eta^{\textbf{e}_{\gamma},\textbf{e}_{Jg_{i}^{c}}}\textbf{e}_{Jg_{i}^{c}}.
\]
To find the value of
$\br{\textbf{e}_{Jg_{i}^{c-1}},\textbf{e}_{Jg_{i}},\textbf{e}_{\gamma}}$
we must compute the degrees of its line bundles. From Section
\ref{review}, we see that
\[
l_{1}=q_{1}-(q_{1}+(c-1)\theta_{1}^{(2)}+q_{1}+\theta_{1}^{(2)}+q_{1}+(a_{2}-1)\theta_{1}^{(1)}+(a_{1}-c-1)\theta_{1}^{(2)})=-1.
\]
\[
l_{2}=q_{2}-(q_{2}+(c-1)\theta_{2}^{(2)}+q_{2}+\theta_{2}^{(2)}+q_{2}+(a_{2}-1)\theta_{2}^{(1)}+(a_{1}-c-1)\theta_{2}^{(2)})=-1.
\]

Therefore, by Axioms 3 and 4 we have that
$\br{\textbf{e}_{Jg_{i}^{c-1}},\textbf{e}_{Jg_{i}},\textbf{e}_{\gamma}}=1$.
In a similar way it can be shown that
$\eta^{\textbf{e}_{\gamma},\textbf{e}_{Jg_{i}^{c}}}=1$, and thus
we find that $(\textbf{e}_{Jg_{i}})^{c}=\textbf{e}_{Jg_{i}^c}$,
which concludes the proof of Lemma 2.8.
\end{proof}

\begin{lemma}
In $(\mathscr{H}_{W,G_{W}},\star)$
\begin{equation}
\textbf{e}_{h_{2}}^{a_{1}}+a_{2}\textbf{e}_{h_{2}}\star \,
  \,\textbf{e}_{h_{1}}^{a_{2}-1}=0\quad \text{and}\quad
\textbf{e}_{h_{1}}^{a_{2}}+a_{1}\textbf{e}_{h_{1}}\star\,
 \,  \textbf{e}_{h_{2}}^{a_{1}-1}=0,
\end{equation}

where $h_{i}=Jg_{i}$, $i \in \{1,2\}$.
\end{lemma}
\begin{proof}

Using the Lemma 2.8, we see that
$\textbf{e}_{h_{2}}^{a_{1}-2}=\textbf{e}_{h_{2}^{a_{1}-2}}$, and
thus
\[
\textbf{e}_{h_{2}}^{a_{1}-1}=\sum_{\alpha,\beta}\br{\textbf{e}_{Jg_{2}^{a_{1}-2}},\textbf{e}_{Jg_{2}},\alpha}\eta^{\alpha,\beta}\beta,
\]
where we need
$\deg_{W}(\textbf{e}_{Jg_{2}^{a_{1}-2}})+\deg_{W}(\textbf{e}_{Jg_{2}})+\deg_{W}(\alpha)=2\hat{c}$.
Making use of \emph{Corollary 2.7}, we note that this last
relation is equivalent to $\deg_W{\alpha}=\hat{c}$. However, this
is only possible if $\alpha \in \mathscr{H}_{i_{W}}$. Therefore,
we have that
\begin{equation}
\textbf{e}_{h_{2}}^{a_{1}-1}=\sum_{\alpha,\beta}
\br{\textbf{e}_{Jg_{2}^{a_{1}-2}},\textbf{e}_{Jg_{2}},\alpha}\eta^{\alpha,\beta}\beta,
\quad \text{where} \quad {\alpha,\beta} \in
\{\left(X_{1}^{a_{1}-1}\right)\textbf{e}_{i_{W}},
\left(X_{2}^{a_{2}-1}\right)\textbf{e}_{i_{W}}\}
\end{equation}
\[
\Rightarrow \textbf{e}_{h_{2}}^{a_{1}}=\sum_{\alpha,\beta}
\br{\textbf{e}_{Jg_{2}^{a_{1}-2}},\textbf{e}_{Jg_{2}},\alpha}\eta^{\alpha,\beta}\br{\beta,\textbf{e}_{Jg_{2}},\textbf{e}_{Jg_{2}^{a_{1}-2}}}\eta^{\textbf{e}_{Jg_{2}^{a_{1}-2}},\delta}\delta
\]
\[
=\Lambda_{0,4}^W(\textbf{e}_{Jg_{2}^{a_{1}-2}},\textbf{e}_{Jg_{2}},\textbf{e}_{Jg_{2}},\textbf{e}_{Jg_{2}^{a_{1}-2}})\eta^{\textbf{e}_{Jg_{2}^{a_{1}-2}},\delta}\delta.
\]
The inverse of $Jg_{2}^{a_{1}-2}$ is given by
$Jg_{1}^{a_{2}-1}g_{2}$, and so
$\delta=\textbf{e}_{Jg_{1}^{a_{2}-1}g_{2}}$. It is not hard to
show that $\eta^{\textbf{e}_{Jg_{2}^{a_{1}-2}},\delta}=1$, which
gives us that
\[
\textbf{e}_{h_{2}}^{a_{1}}=\Lambda_{0,4}^W(\textbf{e}_{Jg_{2}^{a_{1}-2}},\textbf{e}_{Jg_{2}},\textbf{e}_{Jg_{2}},\textbf{e}_{Jg_{2}^{a_{1}-2}})\textbf{e}_{Jg_{1}^{a_{2}-1}g_{2}}.
\]
To find the value of this four-point class, we compute the degrees
of its line bundles,
\[
l_{1}=2q_{1}-(q_{1}+(a_{1}-2)\theta_{1}^{(2)}+q_{1}+\theta_{1}^{(2)}+q_{1}+\theta_{1}^{(2)}+q_{1}+(a_{1}-2)\theta_{1}^{(2)})=-2,
\]
\[
l_{2}=2q_{2}-(q_{2}+(a_{1}-2)\theta_{2}^{(2)}+q_{2}+\theta_{2}^{(2)}+q_{2}+\theta_{2}^{(2)}+q_{2}+(a_{1}-2)\theta_{2}^{(2)})=0.
\]

Using Axiom 5, we find that
$\Lambda_{0,4}^W(\textbf{e}_{Jg_{2}^{a_{1}-2}},\textbf{e}_{Jg_{2}},\textbf{e}_{Jg_{2}},\textbf{e}_{Jg_{2}^{a_{1}-2}})=-a_{2}$,
and thus
\begin{equation}
\textbf{e}_{h_{2}}^{a_{1}}=-a_{2}\textbf{e}_{Jg_{1}^{a_{2}-1}g_{2}}.
\end{equation}

In the same way (25) was obtained, one can show that
\[
\textbf{e}_{h_{1}}^{a_{2}-1}=\sum_{\alpha,\beta}
\br{\textbf{e}_{Jg_{1}^{a_{2}-2}},\textbf{e}_{Jg_{1}},\alpha}\eta^{\alpha,\beta}\beta,
\quad \text{where} \quad {\alpha, \beta} \in
\{\left(X_{1}^{a_{1}-1}\right)\textbf{e}_{i_{W}},
\left(X_{2}^{a_{2}-1}\right)\textbf{e}_{i_{W}}\}
\]
\[
\Rightarrow \textbf{e}_{h_{2}} \star\,
    \, \textbf{e}_{h_{1}}^{a_{2}-1}=\sum_{\alpha,\beta} \br{\textbf{e}_{Jg_{1}^{a_{2}-2}},\textbf{e}_{Jg_{1}},\alpha}\eta^{\alpha,\beta}
\br{\beta,\textbf{e}_{Jg_{2}},\textbf{e}_{Jg_{2}^{a_{1}-2}}}\eta^{\textbf{e}_{Jg_{2}^{a_{1}-2}},\textbf{e}_{Jg_{1}^{a_{2}-1}g_{2}}}\textbf{e}_{Jg_{1}^{a_{2}-1}g_{2}}
\]
\[
=\Lambda_{0,4}^W(\textbf{e}_{Jg_{1}^{a_{2}-2}},\textbf{e}_{Jg_{1}},\textbf{e}_{Jg_{2}},\textbf{e}_{Jg_{2}^{a_{1}-2}})\eta^{\textbf{e}_{Jg_{2}^{a_{1}-2}},\textbf{e}_{Jg_{1}^{a_{2}-1}g_{2}}}\textbf{e}_{Jg_{1}^{a_{2}-1}g_{2}}.
\]
It is not hard to show that
$\eta^{\textbf{e}_{Jg_{2}^{a_{1}-2}},\textbf{e}_{Jg_{1}^{a_{2}-1}g_{2}}}=1$.
Now, to find the value of the four-point class we compute its line
bundle degrees,
\[
l_{1}=2q_{1}-(q_{1}+(a_{2}-2)\theta_{1}^{(1)}+q_{1}+\theta_{1}^{(1)}+q_{1}+\theta_{1}^{(2)}+q_{1}+(a_{1}-2)\theta_{1}^{(2)})=-1,
\]

\[
l_{2}=2q_{2}-(q_{2}+(a_{2}-2)\theta_{2}^{(1)}+q_{2}+\theta_{2}^{(1)}+q_{2}+\theta_{2}^{(2)}+q_{2}+(a_{1}-2)\theta_{2}^{(2)})=-1,
\]

and so by Axioms 3 and 4 we have that
$\Lambda_{0,4}^W(\textbf{e}_{Jg_{1}^{a_{2}-2}},\textbf{e}_{Jg_{1}},\textbf{e}_{Jg_{2}},\textbf{e}_{Jg_{2}^{a_{1}-2}})=1$.
Therefore, we have that
\[
\textbf{e}_{h_{2}} \star \,
  \,
  \textbf{e}_{h_{1}}^{a_{2}-1}=\textbf{e}_{Jg_{1}^{a_{2}-1}g_{2}}.
\]
Putting this together with (26) allows us to show that
$\textbf{e}_{h_{2}}^{a_{1}}+a_{2}\textbf{e}_{h_{2}}\textbf{e}_{h_{1}}^{a_{2}-1}=0$.
Following the steps that led us to this relation, one can show
that
$\textbf{e}_{h_{1}}^{a_{2}}+a_{1}\textbf{e}_{h_{1}}\textbf{e}_{h_{2}}^{a_{1}-1}=0$.
\end{proof}

We are now in a position to prove Theorem 2.5,

\emph{Proof of  Theorem 2.5:}

Consider the map $\varphi :\Q_{W^T}\longrightarrow
\mathscr{H}_{W,G_{W}}$ given by
\[
\overline{X}_{1}^{\alpha}\overline{X}_{2}^{\beta}\longmapsto
\textbf{e}_{Jg_{1}^{\alpha}g_{2}^{\beta}}, \quad
\overline{X}_{1}^{a_{2}-1}\longmapsto
\left(X_{2}^{a_{2}-1}\right)\textbf{e}_{i_{W}},\quad
\overline{X}_{2}^{a_{1}-1}\longmapsto
\left(X_{1}^{a_{1}-1}\right)\textbf{e}_{i_{W}},\quad
\]
where $0 \leq \alpha \leq a_{2}-1$ and $0 \leq \beta \leq
a_{1}-1$.

From Lemma 2.6, it is easy to see that this map is surjective, and
because the dimensions of $\Q_{W^T}$ and $\mathscr{H}_{W,G_{W}}$
are equal, $\varphi$ must be bijective.

Note that $\textbf{e}_{Jg_{i}} \longmapsto \overline{X}_{i}$, $i
\in \{1,2\}$, and that the relations in $\mathscr{Q}_{W^{T}}$ are
given by its Jacobian ideal, i.e.
\begin{equation}
\overline{X}_{2}^{a_{1}}+a_{2}\overline{X}_{2}\overline{X}_{1}^{a_{2}-1}=0,
\quad
\overline{X}_{1}^{a_{2}}+a_{1}\overline{X}_{1}\overline{X}_{2}^{a_{1}-1}=0.
\end{equation}
Therefore, by Lemma 2.9, $\text{Jac}(W^{T}) \subseteq
\text{ker}(\varphi)$, and we have that $\varphi$ is the desired
degree preserving isomorphism.


\bibliographystyle{amsplain}

\providecommand{\bysame}{\leavevmode\hbox
to3em{\hrulefill}\thinspace}

\end{document}